\pdfoutput=1
\RequirePackage{ifpdf}
\ifpdf 
\documentclass[pdftex]{sigma}
\else
\documentclass{sigma}
\fi

\usepackage[ansinew]{inputenc}
\newtheorem{Theorem}{Theorem}[section]

\newtheorem{Lemma}[Theorem]{Lemma}
\newtheorem{Proposition}[Theorem]{Proposition}
{ \theoremstyle{definition}
\newtheorem{Definition}[Theorem]{Definition}

 }
\numberwithin{equation}{section}

\DeclareMathOperator{\supp}{supp}

\begin{document}

\newcommand{\abs}[1]{\left\vert #1 \right\vert}

\newcommand{\set}[1]{\left\{ #1 \right\}}

\newcommand{\seq}[1]{\langle #1 \rangle}

\newcommand{\arXivNumber}{1411.2000}

\allowdisplaybreaks

\renewcommand{\thefootnote}{$\star$}

\renewcommand{\PaperNumber}{026}

\FirstPageHeading

\ShortArticleName{On the $q$-Charlier Multiple Orthogonal Polynomials}

\ArticleName{On the $\boldsymbol{q}$-Charlier Multiple Orthogonal Polynomials\footnote{This paper is a~contribution to
the Special Issue on Exact Solvability and Symmetry Avatars in honour of Luc Vinet.
The full collection is available at
\href{http://www.emis.de/journals/SIGMA/ESSA2014.html}{http://www.emis.de/journals/SIGMA/ESSA2014.html}}}

\Author{Jorge ARVES\'U and Andys M.~RAM\'{I}REZ-ABERASTURIS}
\AuthorNameForHeading{J.~Arves\'u and A.M.~Ram\'{\i}rez-Aberasturis}

\Address{Department of Mathematics, Universidad Carlos III de Madrid,\\
Avenida de la Universidad, 30, 28911, Legan\'es, Spain}
\Email{\href{mailto:jarvesu@math.uc3m.es}{jarvesu@math.uc3m.es}, \href{mailto:aramirez@math.uc3m.es}{aramirez@math.uc3m.es}}
\URLaddress{\url{http://gama.uc3m.es/index.php/jarvesu.html}}

\ArticleDates{Received November 10, 2014, in f\/inal form March 23, 2015; Published online March 28, 2015}

\Abstract{We introduce a~new family of special functions, namely $q$-Charlier multiple orthogonal polynomials.
These polynomials are orthogonal with respect to $q$-analogues of Poisson distributions.
We focus our attention on their structural properties.
Raising and lowering operators as well as Rodrigues-type formulas are obtained.
An explicit representation in terms of a $q$-analogue of the second of Appell's hypergeometric functions is given.
A~high-order linear $q$-dif\/ference equation with polynomial coef\/f\/icients is deduced.
Moreover, we show how to obtain the nearest neighbor recurrence relation from some dif\/ference operators involved in the
Rodrigues-type formula.}

\Keywords{multiple orthogonal polynomials; Hermite--Pad\'e approximation; dif\/ference equations; classical orthogonal
polynomials of a~discrete variable; Charlier polynomials; $q$-poly\-no\-mials}

\Classification{42C05; 33E30; 33C47; 33C65}

\renewcommand{\thefootnote}{\arabic{footnote}} \setcounter{footnote}{0}

\section{Introduction}

Let $\vec{\mu}=(\mu_{1},\ldots,\mu_{r})$ be a~vector of~$r$ positive Borel measures on $\mathbb{R}$, and let
$\vec{n}=(n_{1},\ldots,n_{r})\in\mathbb{N}^{r}$ be a~multi-index.
By $\mathbb{N}$ we denote the set of all nonnegative integers.
A~type~II multiple orthogonal polynomial $P_{\vec{n}}$, corresponding to the multi-index $\vec{n}$, is a~polynomial of
degree $\leq\vert \vec{n}\vert =n_{1}+\dots +n_{r}$ which satisf\/ies the orthogonality
conditions~\cite{Nikishin}
\begin{gather}
\int_{\Omega_{i}}P_{\vec{n}}(x) x^{k}d\mu_{i}(x) =0,
\qquad
k=0,\ldots,n_{i}-1,
\qquad
i=1,\ldots,r,
\label{OrthC}
\end{gather}
where $\Omega_{i}$ is the smallest interval that contains $\supp(\mu_{i}) $.
In this paper we will consider the situation when $P_{\vec{n}}$ is a~monic multiple orthogonal polynomial and has
exactly degree $|\vec{n}|$.
If the measures in~\eqref{OrthC} are discrete
\begin{gather}
\mu_{i}=\sum\limits_{k=0}^{N_{i}}\omega_{i,k}\delta_{x_{i,k}},
\qquad
\omega_{i,k}>0,
\qquad
x_{i,k}\in \mathbb{R},
\qquad
N_{i}\in \mathbb{N\cup}\{+\infty \},
\qquad
i=1,2,\ldots,r,
\label{dismeasure}
\end{gather}
where $\delta_{x_{i,k}}$ denotes the Dirac delta function and $x_{i_{1},k}\neq x_{i_{2},k}$, $k=0,\ldots,N_{i}$,
whenever $i_{1}\neq i_{2}$, the corresponding polynomial solution is called discrete multiple orthogonal polynomial.

In~\cite{arvesu_vanAssche} was studied some type II discrete multiple orthogonal polynomials on the linear lattice
$x(s)=s$.
In particular, multiple Charlier polynomials were considered (when the component measures of $\vec{\mu}$ are dif\/ferent
Poisson distributions).
In this paper we will introduce a $q$-analogue of such multiple orthogonal polynomials (when the component measures of
$\vec{\mu}$ are dif\/ferent $q$-Poisson distributions) and study their algebraic properties (structural properties).
We are motivated by the recent applications that have been found for their predecessors: Multiple Charlier polynomials.
In particular, these polynomials appear in remainder Pad\'e approximation for the exponential
function~\cite{Prevost-Rivoal} and as common eigenstate of a~set of~$r$ non-Hermitian oscillator
Hamiltonians~\cite{miki-vinet-zhedanov}.
Furthermore, in~\cite{Ndayiragije-Assche} the authors pointed out a~possible relationship of these polynomials to the
orthogonal functions appearing in two speed totally asymmetric simple exclusion process (TASEP)~\cite{borodin-2}.
Our new $q$-family of multiple orthogonal polynomials is likely to be of relevance in TASEP-like models as well as
in $q$-extensions of the mathematical problems addressed in~\cite{miki-vinet-zhedanov,Prevost-Rivoal}.

Recently, two families of $q$-multiple orthogonal polynomials and some of their structural properties have been
studied~\cite{arvesu-qHahn,Postelmans_Assche}.
Moreover, in~\cite{arvesu-esposito,lee,Assche_diff-eq} an $(r+1)$-order dif\/ference equation for some discrete multiple
orthogonal polynomials was obtained.
An interesting fact is that these polynomials are common eigenfunctions of two distinct linear dif\/ference operators of
order $(r+1)$ since multiple orthogonal polynomials also satisfy $(r+2)$-term recurrence relations.
The explicit expressions for the coef\/f\/icients of this relation are the main ingredient for the study of some type of
asymptotic behaviors for these polynomials.
For instance, in~\cite{apt-arv} the weak asymptotics was studied for multiple Meixner polynomials of the f\/irst and
second kind, respectively.
The zero distribution of multiple Meixner polynomials was also studied.
In~\cite{Ndayiragije-Assche} the recurrence relation is also used for studying the ratio asymptotic behavior of multiple
Charlier polynomials (introduced in~\cite{arvesu_vanAssche}) and from it the authors obtain the asymptotic distribution
of the zeros after a~suitable rescaling.
Furthermore, this recurrence relation is a~key-ingredient for attaining a~Christof\/fel--Darboux
kernel~\cite{daems_kuijlaars_kernel} among other applications, which plays important role in correlation kernel (for
instance in the unitary random matrix model with external source).

In this paper, we will mainly focus on two structural properties satisf\/ied by the multiple orthogonal polynomials
introduced here, namely recurrence relations and the dif\/ference equation (with respect to the independent variable).
The asymptotic analysis is out of the scope of this paper.
Although we plan to address this question in a~future publication.

\looseness=-1
The structure of the paper is as follows.
Section~\ref{sec1} introduces the context and the background materials.
In Section~\ref{q-ch} we will def\/ine the $q$-Charlier multiple orthogonal polynomials (new in the literature) and obtain
raising operators and then the Rodrigues-type formula for these special functions when~$r$ orthogonality conditions are
considered.
The explicit series representation in terms of a $q$-analogue of the second of Appell's hypergeometric functions is
given for $\vec{n}=(n_1,n_2)$.
In Section~\ref{df-eq}, an $(r+1)$-order $q$-dif\/ference equation on a~non-uniform lattice $x(s)$ is obtained.
Section~\ref{rr} deals with the $(r+2)$-term recurrence relations.
In particular, the nearest neighbor recurrence relation is obtained.
Explicit expressions for the recurrence coef\/f\/icients are given.
Finally, in Section~\ref{conclu} some of our f\/indings are summarized.
A~special limiting case when the parameter~$q$ involved in the components of the vector measure tends to one is pointed out.
Under this limit the corresponding structural properties for the multiple Charlier polynomials are recovered.

\section{Multiple Charlier polynomials}\label{sec1}

For the discrete measures~\eqref{dismeasure} we have that $\supp(\mu_{i})$ is the closure of
$\{x_{i,k}\}_{k=0}^{N_{i}}$ and that $\Omega_{i}$ is the smallest closed interval on $\mathbb{R}$ which
contains $\{x_{i,k}\}_{k=0}^{N_{i}}$.
Moreover, the above orthogonality conditions~\eqref{OrthC} give a~linear system of $\vert\vec{n}\vert$
homogeneous equations for the $\vert\vec{n}\vert +1$ unknown coef\/f\/icients of $P_{\vec{n}}(x)$.
This polynomial solution $P_{\vec{n}}$ always exists.
We focus our attention on a~unique solution (up to a~multiplicative factor) with $\deg P_{\vec{n}}(x)=|\vec{n}|$.
If this happen for every multi-index $\vec{n}$, we say that $\vec{n}$ is normal~\cite{Nikishin}.
If the above system of measures forms an $AT$ system~\cite{Nikishin} then every multi-index is normal.
Indeed, we will deal with such system of discrete measures, where $\Omega_{i}=\Omega$ for each $i=1,2,\ldots,r$.

\begin{Definition}
The system of positive discrete measures $\mu_{1},\mu_{2},\ldots,\mu_{r}$, given in~\eqref{dismeasure}, forms an $AT$
system if there exist~$r$ continuous functions $\upsilon_{1},\ldots,\upsilon_{r}$ on~$\Omega$ with
$\upsilon_{i}(x_{k})=\omega_{i,k}$, $k=0,\ldots,N_{i}$, $i=1,2,\ldots,r$, such that the
$\vert\vec{n}\vert$ functions
\begin{gather*}
\upsilon_{1}(x),x\upsilon_{1}(x),\ldots,x^{n_{1}-1}\upsilon_{1}(x),\ldots,\upsilon_{r}(x),x\upsilon_{r}(x),\ldots,x^{n_{r}-1}\upsilon_{r}(x),
\end{gather*}
form a~Chebyshev system on~$\Omega$ for each multi-index $\vec{n}$ with $\vert \vec{n}\vert <N+1$, i.e.,
every linear combination $\sum\limits_{i=1}^{r}Q_{n_{i}-1}(x) \upsilon_{i}(x)$, where $Q_{n_{i}-1}\in
\mathbb{P}_{n_{i}-1}\setminus \{0\}$, has at most $\vert \vec{n}\vert -1$ zeros on~$\Omega$.
\end{Definition}

The monic multiple Charlier polynomials~\cite{arvesu_vanAssche} $C_{\vec{n}}^{\vec{\alpha}}(x)$, with multi-index
$\vec{n}\in \mathbb{N}^{r}$ and deg\-ree~$\vert \vec{n}\vert $ satisfy the following orthogonality conditions
with respect to~$r$ Poisson distributions with dif\/ferent positive parameters $\alpha_{1},\ldots,\alpha_{r}$ (indexed~by
$\vec{\alpha}=(\alpha_{1},\ldots,\alpha_{r})$)
\begin{gather*}
\sum\limits_{x=0}^{\infty}C_{\vec{n}}^{\vec{\alpha}}(x)(-x)_{j}
\upsilon^{\alpha_{i}}(x)=0,
\qquad
j=0,\ldots,n_{i}-1,
\qquad
i=1,\ldots,r,
\end{gather*}
where
\begin{gather*}
\upsilon^{\alpha_{i}}(x)=
\begin{cases}
\dfrac{\alpha_{i}^{x}}{\Gamma (x+1)}, & \text{if} \quad  x\in\mathbb{R}\setminus\mathbb{Z}_{-},
\\
0, & \text{otherwise,}
\end{cases}
\end{gather*}
and $(x)_j=(x)(x+1)\cdots(x+j-1)$, $(x)_0=1$, $j\geq 1$, denotes
the Pochhammer symbol.

In~\cite{arvesu_vanAssche} the authors consider a~normal multi-index $\vec{n}\in\mathbb{N}^{r}$, whenever
$\alpha_{i}>0$, $i=1,2,\ldots,r$, and with all the $\alpha_{i}$ dif\/ferent.
Moreover, it was found the following raising operators
\begin{gather}
\mathcal{L}_{\vec{n}}^{\alpha_{i}}\big[C_{\vec{n}}^{\vec{\alpha}}(x)\big ]
=-C_{\vec{n}+\vec{e}_{i}}^{\vec{\alpha}}(x),
\qquad
i=1,\ldots,r,
\label{raisingx}
\end{gather}
where
\begin{gather*}
\mathcal{L}_{\vec{n}}^{\alpha_{i}}\overset{\rm def}{=}\frac{\alpha_{i}}{\upsilon^{\alpha_{i}}(x)}\bigtriangledown\upsilon^{\alpha_{i}}(x),
\end{gather*}
and $\bigtriangledown f(x)=f(x)-f(x-1)$ denotes the backward dif\/ference operator.
As a~consequence of~\eqref{raisingx} there holds the Rodrigues-type formula
\begin{gather}
C_{\vec{n}}^{\vec{\alpha}}(x)
=(-1)^{\vert\vec{n}\vert}\left(\prod\limits_{j=1}^{r}\alpha_{j}^{n_{j}}\right) \Gamma
(x+1)\mathcal{C}_{\vec{n}}^{\vec{\alpha}}\left(\frac{1}{\Gamma(x+1)}\right),
\label{Rodrigues-multi}
\end{gather}
where
\begin{gather*}
\mathcal{C}_{\vec{n}}^{\vec{\alpha}}=\prod\limits_{i=1}^{r}
\big(\alpha_{i}^{-x}\bigtriangledown^{n_{i}}\alpha_{i}^{x}\big).
\end{gather*}

Two important structural properties are known for multiple Charlier polynomials~\cite{arvesu_vanAssche}, namely the
$(r+1)$-order linear dif\/ference equation~\cite{lee}
\begin{gather}
\prod\limits_{i=1}^{r}\mathcal{L}_{\vec{n}}^{\alpha_{i}}\big[\bigtriangleup C_{\vec{n}}^{\vec{\alpha}}(x) \big]
+\sum\limits_{i=1}^{r}n_{i}\prod\limits_{\substack{j=1
\\
j\neq i}}^{r} \mathcal{L}_{\vec{n}}^{\alpha_{j}}\big[C_{\vec{n}}^{\vec{\alpha}}(x)\big] =0,
\label{opdi-1}
\end{gather}
where $\bigtriangleup f(x)=f(x+1)-f(x)$,
and the recurrence relation~\cite{arvesu_vanAssche,Hane}
\begin{gather}
xC_{\vec{n}}^{\vec{\alpha}}(x) =C_{\vec{n}+\vec{e_{k}}}^{\vec{\alpha}}(x) +(\alpha_{k}+\vert
\vec{n}\vert ) C_{\vec{n}}^{\vec{\alpha}}(x)
+\sum\limits_{i=1}^{r}\alpha_{i}n_{i}C_{\vec{n}-\vec{e}_{i}}^{\vec{\alpha}}(x),
\label{AR}
\end{gather}
where the multi-index $\vec{e}_{i}$ is the standard~$r$ dimensional unit vector with the~$i$-th entry equals $1$ and $0$
otherwise.

Interestingly, the multiple Charlier polynomials $C_{\vec{n}}^{\vec{\alpha}}(x)$ are common eigenfunctions of the above
two linear dif\/ference operators of order $(r+1)$, namely~\eqref{opdi-1} and~\eqref{AR}.

\section[$q$-Charlier multiple orthogonal polynomials]{$\boldsymbol{q}$-Charlier multiple orthogonal polynomials}\label{q-ch}

Let us begin by recalling the def\/inition of $q$-multiple orthogonal polynomials~\cite{arvesu-qHahn}.
\begin{Definition}
\label{defqpo}
A~polynomial $P_{\vec{n}}(x(s))$ on the lattice $x(s)=c_1q^s+c_3$, $q\in\mathbb{R}^{+}\setminus\{1\}$,
$c_1,c_3\in\mathbb{R}$, is said to be a $q$-multiple orthogonal polynomial of a~multi-index $\vec{n} \in \mathbb{N}^{r}$
with respect to positive discrete measures $\mu_{1},\mu_{2},\ldots,\mu_{r}$
such that $\supp(\mu_{i})\subset \Omega_{i}\subset \mathbb{R}$, $i=1,2,\ldots,r$, if there hold conditions
\begin{gather}
(a)
\quad
\deg P_{\vec{n}}(x(s) ) \leq \vert \vec{n}\vert =n_{1}+n_{2}+\dots +n_{r},
\nonumber
\\
(b)
\quad
\sum\limits_{s=0}^{N_{i}}P_{\vec{n}}(x(s) ) x(s)^{k}d\mu_{i}=0,
\qquad
k=0,\ldots,n_{i}-1,
\qquad
N_{i}\in\mathbb{N\cup}\{+\infty\}.
\label{OrthC2}
\end{gather}
\end{Definition}

Let us consider the following~$r$ positive discrete measures on $\mathbb{R}^{+}$,
\begin{gather}
\mu_{i}=\sum\limits_{s=0}^{\infty}\omega_{i}(k)\delta(k-s),
\qquad
\omega_{i}>0,
\qquad
i=1,2,\ldots,r.
\label{Measu}
\end{gather}
Here $\omega_{i}(s)=\upsilon_{q}^{\alpha_{i}}(s) \bigtriangleup x(s-1/2)$, which involve the $q$-analogue of
Poisson distributions
\begin{gather*}
\upsilon_{q}^{\alpha_{i}}(s) =
\begin{cases}
\dfrac{\alpha_{i}^{s}}{\Gamma_{q}(s+1)}, & \text{if}\quad s\in \mathbb{R}^{+}\cup \{0\},
\\
0, &\text{otherwise},
\end{cases}
\end{gather*}
where $\alpha_{i}>0$, $i=1,2,\ldots,r$, with all the $\alpha_{i}$ dif\/ferent.
Recall that the $q$-Gamma function is given~by
\begin{gather}
\Gamma_q(s) =
\begin{cases}
f(s;q)=(1-q)^{1-s} \dfrac{\prod\limits_{k\geq0} (1-q^{k+1})}{\prod\limits_{k\geq0} (1-q^{s+k})}, & 0<q<1,
\\
q^{\frac{(s-1)(s-2)}{2}}f\big(s;q^{-1}\big), & q>1.
\end{cases}
\label{q-gamma-clas}
\end{gather}
See also~\cite{Gasper,nsu} for the above def\/inition of the $q$-Gamma function.

\begin{Lemma}
\label{le-1}
The system of functions
\begin{gather}
 \alpha_1^s, x(s)\alpha_1^s,\ldots,x(s)^{n_1-1}\alpha_1^s,\ldots, \alpha_r^s,
x(s)\alpha_r^s,\ldots,x(s)^{n_r-1}\alpha_r^s,
\label{cheb}
\end{gather}
with $\alpha_{i}>0$, $i=1,2,\ldots,r$, and $(\alpha_i/\alpha_j)\neq q^k$,
$k\in\mathbb{Z}$,
$i,j=1,\dots,r$,
$i\neq j$, forms a~Chebyshev system on $\mathbb{R}^{+}$ for every
$\vec{n}=(n_1,\dots,n_r)\in\mathbb{N}^{r}$.
\end{Lemma}

\begin{proof}
This means that every linear combination $\sum\limits_{i=1}^{r}Q_{n_{i}-1}(x(s))\alpha_{i}^s$ has at most
$\vert\vec{n}\vert-1$ zeros on $\mathbb{R}^{+}$ for every $Q_{n_{i}-1}(x(s))\in
\mathbb{P}_{n_{i}-1}\setminus \{0\}$.
Since $x(s)=c_1q^s+c_3$, where $c_1$, $c_3$ are constants, we consider
$\sum\limits_{i=1}^{r}Q_{n_{i}-1}(q^s)\alpha_{i}^s $, instead.
Thus, the system~\eqref{cheb} transforms into
\begin{gather*}
a_{1,0}^s, a_{1,1}^s,\ldots,a_{1,n_1-1}^s,\ldots, a_{r,0}^s, a_{r,1}^s,\ldots,a_{r,n_r-1}^s,
\end{gather*}
where $a_{i,k}=(q^k\alpha_i)$, with $k=0,\dots,n_i-1$, $i=1,\dots,r$.
Observe that $a_{j,m}\neq a_{l,p}$ for $j\neq l$, $m\neq p$.
Hence, identity $a_{i,k}=e^{\log a_{i,k}}$ yields the well-known Chebyshev system (see~\cite[p.~138]{Nikishin})
\begin{gather*}
e^{s\log a_{1,0}}, e^{s\log a_{1,1}},\ldots,e^{s\log a_{1,n_1-1}},\ldots, e^{s\log a_{r,0}}, e^{s\log
a_{r,1}},\ldots,e^{s\log a_{r,n_r-1}}.
\end{gather*}
Then, we conclude that the functions~\eqref{cheb} form a~Chebyshev system on $\mathbb{R}^{+}$.
\end{proof}

As a~consequence of Lemma~\ref{le-1} the system of measures $\mu_1,\mu_2,\ldots,\mu_r$ given in~\eqref{Measu} forms an
AT system on $\mathbb{R}^{+}$.
Using this fact, we will rewrite Def\/inition~\ref{defqpo} for this system of measures.
However, we f\/irst def\/ine the $q$-analogue of the Stirling polynomials denoted by $[s]_{q}^{(k)}$, which is
a~polynomial of degree~$k$ in the variable $x(s)=(q^s-1)/(q-1)$, as follows
\begin{gather*}
[s]_{q}^{(k)}=\prod\limits_{j=0}^{k-1}\frac{q^{s-j}-1}{q-1} =x(s) x(s-1) \cdots x(s-k+1)
\qquad\!
\text{for}
\quad
k>0,
\qquad\!
\text{and}
\qquad\!
[s]_{q}^{(0)}=1.
\end{gather*}

\begin{Definition}
A~polynomial $C_{q,\vec{n}}^{\vec{\alpha}}(s)$, with multi-index $\vec{n}\in \mathbb{N}^{r}$ and degree $\vert
\vec{n}\vert$ that verif\/ies the orthogonality conditions
\begin{gather}
\sum\limits_{s=0}^{\infty}C_{q,\vec{n}}^{\vec{\alpha}}(s) [s]_{q}^{(k)}\upsilon_{q}^{\alpha_{i}}(s)
\bigtriangleup x(s-1/2) =0,
\qquad
0\leq k\leq n_{i}-1,
\qquad
i=1,\ldots,r,
\qquad
\label{NOrthC}
\end{gather}
(see~\eqref{OrthC2} with respect to the measures~\eqref{Measu}) is said to be the $q$-Charlier multiple orthogonal
polynomial.
\end{Definition}

Let us point out some observations derived from the above def\/inition.
When $r=1$ we recover the scalar $q$-Charlier polynomials~\cite{alv_arv}.
The orthogonality conditions~\eqref{OrthC2} have been written more conveniently as~\eqref{NOrthC} since
$[s]_q^{(k)}=q^{-\binom{k}{2}}x^k(s)+\text{lower terms}=\mathcal{O}(q^{ks})$.
Here, the symbol $\mathcal{O}(\cdot)$ stands for big-O notation.
Indeed, when~$q$ goes to 1, the symbol $[s]_{q}^{(k)}$ converges to $(-1)^{k}(-s)_{k}$, where $(s)_{k}$
denotes the Pochhammer symbol.
Hence, one can recover the multiple Charlier polynomials given in~\cite{arvesu_vanAssche} as a~limiting case, provided
that the lattice $x(s)=(q^s-1)/(q-1)$.
Moreover, in the sequel we will only use this lattice.
Finally, we have an AT-system of positive discrete measures, then the $q$-Charlier multiple orthogonal polynomial
$C_{q,\vec{n}}^{\vec{\alpha}}(s)$ has exactly $\vert \vec{n}\vert$ dif\/ferent zeros on $\mathbb{R}^{+}$
(see~\cite[Theorem 2.1, pp.~26--27]{arvesu_vanAssche}).

For monic $q$-Charlier multiple orthogonal polynomials we have~$r$ raising operators
\begin{gather}
\mathcal{D}_{q}^{\alpha_{i}}C_{q,\vec{n}}^{\vec{\alpha}}(s)
=-q^{1/2}C_{q,\vec{n}+\vec{e}_{i}}^{\vec{\alpha}_{i,1/q}}(s),
\qquad
i=1,\ldots,r,
\label{ROpqMCharlier}
\end{gather}
where
\begin{gather*}
\mathcal{D}_{q}^{\alpha_{i}}\overset{\rm def}{=}
\left(\frac{\alpha_{i}q^{|\vec{n}|}}{\upsilon_{q}^{\alpha_{i}/q}(s)} \nabla \upsilon_{q}^{\alpha_{i}}(s)\right),
\qquad
\nabla \overset{\rm def}{=} \frac{\bigtriangledown}{\bigtriangledown x(s+1/2)},
\qquad
\vec{\alpha}_{i,1/q}=(\alpha_{1},\ldots,\alpha_{i}/q,\ldots,\alpha_{r}).
\end{gather*}
Furthermore,
\begin{gather*}
q^{-\vert \vec{n}\vert -1/2}\mathcal{D}_{q}^{\alpha_{i}}f(s) =[\alpha_{i}-x(s) ] f(s) +x(s)
\bigtriangledown f(s),
\end{gather*}
for any function $f(s)$ def\/ined on the discrete variable~$s$.
Notice that we call $\mathcal{D}_{q}^{\alpha_{i}}$ a~raising operator since the~$i$-th component of the multi-index
$\vec{n}$ in~\eqref{ROpqMCharlier} is increased by 1.
For f\/inding this operator we have replace $[s]_{q}^{(k)}$ in~\eqref{NOrthC} by the following f\/inite-dif\/ference
expression
\begin{gather}
[s]_{q}^{(k)}=\frac{q^{(k-1) /2}}{[k+1]_{q}^{(1)}}\nabla[s+1]_{q}^{(k+1)},
\label{diffs}
\end{gather}
and then used summation by parts along with conditions
$\upsilon_{q}^{\alpha_{i}}(-1)=\upsilon_{q}^{\alpha_{i}}(\infty)=0$.

In the sequel we will only consider monic $q$-Charlier multiple orthogonal polynomials.
In addition, the following dif\/ference operators will appear regularly
\begin{gather}
\Delta  \overset{\rm def}{=}\frac{\bigtriangleup} {\bigtriangleup x(s-1/2)},
\label{lower}
\\
\nabla^{n_{i}} =\underbrace{\nabla\cdots\nabla}_{n_{i}~\text{times}},
\label{nbackward}
\end{gather}
and $\bigtriangledown x_{1}(s)\overset{\rm def}{=}\bigtriangledown x(s+1/2)= \bigtriangleup x(s-1/2)=q^{s-1/2}$.

\begin{Proposition}
There holds the following $q$-analogue of Rodrigues-type formula
\begin{gather}
C_{q,\vec{n}}^{\vec{\alpha}}(s)
=\mathcal{K}_{q}^{\vec{n},\vec{\alpha}}\Gamma_{q}(s+1)\mathcal{C}_{q,\vec{n}}^{\vec{\alpha}}\left(\frac{1}{\Gamma_{q}(s+1)}\right),
\label{RFormula}
\end{gather}
where
\begin{gather}
\mathcal{C}_{q,\vec{n}}^{\vec{\alpha}}=\prod\limits_{i=1}^{r}\mathcal{C}_{q,n_{i}}^{\alpha_{i}},
\qquad
\mathcal{C}_{q,n_{i}}^{\alpha_{i}}=(\alpha_{i})^{-s} \nabla^{n_{i}}(\alpha_{i}q^{n_i})^{s},
\label{Op-q-Ch}
\end{gather}
and
\begin{gather*}
\mathcal{K}_{q}^{\vec{n},\vec{\alpha}}=(-1)^{\vert \vec{n} \vert}q^{-\frac{\vert \vec{n}
\vert}{2}}\left(\prod\limits_{i=1}^{r}\alpha_i^{n_i}\right)\left(\prod\limits_{i=1}^{r}q^{n_i\sum\limits_{j=i}^{r}n_{j}}\right).
\end{gather*}

\end{Proposition}

\begin{proof}
Replacing $[s]_{q}^{(k)}$ in~\eqref{NOrthC} by the right-hand side of expression~\eqref{diffs} one has
\begin{gather*}
\sum\limits_{s=0}^{\infty}C_{q,\vec{n}}^{\vec{\alpha}}(s) \nabla
[s+1]_{q}^{(k+1)}\upsilon_{q}^{\alpha_{i}}(s) \bigtriangledown x_{1}(s) =0,
\qquad
0\leq k\leq n_{i}-1,
\qquad
i=1,\ldots,r.
\qquad
\end{gather*}
Hence, using summation by parts and the raising operators~\eqref{ROpqMCharlier} in a~recursive way one obtains the
Rodrigues-type formula~\eqref{RFormula}.
\end{proof}

Notice that equation~\eqref{RFormula} provides an explicit expression for the monic $q$-Charlier multiple orthogonal
polynomials.
Indeed, by using formula~(3.2.29) from~\cite{nsu} as follows
\begin{gather*}
\nabla^{m}f(s)=q^{\binom{m+1}{2}/2-ms} \sum\limits_{k=0}^{m}{m\brack k}(-1)^k q^{\binom{m-k}{2}}f(s-k),
\\
{m\brack k}=\frac{(q;q)_{m}}{(q;q)_{k}(q;q)_{m-k}},
\qquad
m=1,2,\ldots,
\end{gather*}
where
\begin{gather*}
(a;q)_k=\prod\limits_{j=0}^{k-1}(1-aq^{j})
\qquad
\text{for}
\quad
k>0,
\qquad
\text{and}
\qquad
(a;q)_0=1,
\end{gather*}
denotes the $q$-analogue of the Pochhammer symbol~\cite{Gasper, Koek, nsu}, one obtains the following relation for
the $q$-Charlier multiple orthogonal polynomials (for multi-index $\vec{n}=(n_1, n_2)$)
\begin{gather}
C_{q,n_1,n_2}^{\alpha_1,\alpha_2}(s)=
(-\alpha_1)^{n_1}(-\alpha_2)^{n_2}q^{n_1^2+n_2^2+n_1n_2-(n_1+n_2)/2}
\nonumber
\\
\hphantom{C_{q,n_1,n_2}^{\alpha_1,\alpha_2}(s)=}{}
 \times \left(\frac{\Gamma_q(s+1)}{\alpha_2^s}\nabla^{n_2}(\alpha_2q^{n_2})^s\right)
\left(\frac{1}{\alpha_1^s}\nabla^{n_1}\frac{(\alpha_1q^{n_1})^s}{\Gamma_q(s+1)}\right)
\nonumber
\\
\hphantom{C_{q,n_1,n_2}^{\alpha_1,\alpha_2}(s)}{}
=
(-\alpha_1)^{n_1}(-\alpha_2)^{n_2}q^{n_1^2+n_2^2+n_1n_2-(n_1+n_2-\binom{n_{1}+1}{2}-\binom{n_{2}+1}{2})/2}
\nonumber
\\
\hphantom{C_{q,n_1,n_2}^{\alpha_1,\alpha_2}(s)=}{}
\times
   \sum\limits_{k=0}^{n_1}\sum\limits_{l=0}^{n_2}(-1)^{l+k} {n_2\brack l}{n_1\brack
k}\frac{q^{\binom{n_2-l}{2}-ln_{2}+\binom{n_{1}-k}{2}-kn_{1}}}{\alpha_{2}^{l}\alpha_{1}^{k}}
\frac{\Gamma_q(s+1)}{\Gamma_q(s-k-l+1)}
\label{eqsph2}
.
\end{gather}
Observe that from~\eqref{q-gamma-clas} we have
\begin{gather*}
\frac{\Gamma_{q}(s+1)}{\Gamma_{q}(s-l-k+1)}=[s]^{(k+l)}_{q}=\frac{(q^{-s};q)_{l+k}}{(q-1)^{l+k}}q^{s(k+l)-\binom{k+l}{2}}.
\end{gather*}
Now, using relation
\begin{gather*}
(q^{-n};q)_{k}=\frac{(q;q)_{n}}{(q;q)_{n-k}}(-1)^{k}q^{\binom{k}{2}-nk},
\qquad
k=0,1,\ldots,
\end{gather*}
(see \cite[formula~(1.8.18)]{Koek}) the above expression~\eqref{eqsph2} transforms into
\begin{gather}
C_{q,n_1,n_2}^{\alpha_1,\alpha_2}(s) =
(-\alpha_1)^{n_1}(-\alpha_2)^{n_2}q^{n_1^2+n_2^2+n_1n_2+\frac{3}{2}\left(\binom{n_{1}}{2}+\binom{n_{2}}{2}\right)}
\nonumber
\\
\hphantom{C_{q,n_1,n_2}^{\alpha_1,\alpha_2}(s) =}{}\times
   \sum\limits_{k=0}^{n_1}\sum\limits_{l=0}^{n_2}
\frac{(q^{-s};q)_{k+l}(q^{-n_1};q)_{k}(q^{-n_2};q)_{l}}{q^{\binom{k+l}{2}}(q;q)_{k}(q;q)_{l}}
\left(\frac{q^{s+1-n_1}}{\alpha_1(q-1)}\right)^{k}\left(\frac{q^{s+1-n_2}}{\alpha_2(q-1)}\right)^{l}.
\label{SHYP}
\end{gather}
Finally, from~\eqref{SHYP} we have
\begin{gather}
C_{q,n_1,n_2}^{\alpha_1,\alpha_2}(s)= (-\alpha_1)^{n_1}(-\alpha_2)^{n_2}
q^{n_1^2+n_2^2+n_1n_2+\frac{3}{2}\left(\binom{n_{1}}{2}+\binom{n_{2}}{2}\right)}
\nonumber
\\
\hphantom{C_{q,n_1,n_2}^{\alpha_1,\alpha_2}(s)=}{}\times
   \lim_{\gamma\rightarrow+\infty} \Phi_{2}\left(q^{-s};q^{-n_1},q^{-n_2};\gamma,\gamma;\frac{\gamma
q^{s+1-n_1}}{\alpha_1(1-q)}, \frac{\gamma q^{s+1-n_2}}{\alpha_2(1-q)}\right),
\label{phi2_charli}
\end{gather}
where
\begin{gather*}
\Phi_{2}(\zeta;\beta,\beta';\gamma,\gamma';x,y)=\sum\limits_{m=0}^{+\infty}
\sum\limits_{n=0}^{+\infty}\frac{(\zeta;q)_{m+n}(\beta;q)_{m}(\beta';q)_{n}}{(\gamma;q)_{m}
(\gamma';q)_{n}(q;q)_{m}(q;q)_{n}}q^{-mn}x^my^n,
\end{gather*}
is a $q$-analogue of the second of Appell's hypergeometric functions of two variables (see \cite[formula~(23), p.~89]{Ernst}).

Alternatively, in~\eqref{SHYP} the $q$-analogue of the Pochhammer symbol can be rewritten in terms of the $q$-falling
factorials, which allows to express the $q$-Charlier multiple orthogonal polynomials in terms of the selected basis
$[s]_q^{(k)}$, $k=0,1,\ldots$, i.e.~
\begin{gather}
C_{q,n_1,n_2}^{\alpha_1,\alpha_2}(s) =
(-\alpha_1)^{n_1}(-\alpha_2)^{n_2}q^{n_1^2+n_2^2+n_1n_2+\frac{3}{2}\left(\binom{n_1}{2}+\binom{n_2}{2}\right)}
\nonumber
\\
\hphantom{C_{q,n_1,n_2}^{\alpha_1,\alpha_2}(s) =}{}\times
   \sum\limits_{k=0}^{n_1}\sum\limits_{l=0}^{n_2}
\frac{[s]_q^{(k+l)}[n_1]_q^{(k)}
[n_2]_q^{(l)}}{[k]_q![l]_q!}
q^{\binom{k+1}{2}+\binom{l+1}{2}}\left(\frac{-1}{q^{2n_1}\alpha_1}\right)^k\left(\frac{-1}{q^{2n_2}\alpha_2}\right)^l.
\label{SHYP2}
\end{gather}
If we replace~$k$ by $j-l$ in~\eqref{SHYP2}, we obtain
\begin{gather}
C_{q,n_1,n_2}^{\alpha_1,\alpha_2}(s) =
 (-\alpha_1)^{n_1}(-\alpha_2)^{n_2}q^{n_1^2+n_2^2+n_1n_2+\frac{3}{2}\left(\binom{n_1}{2}+\binom{n_2}{2}\right)}
\nonumber
\\
\hphantom{C_{q,n_1,n_2}^{\alpha_1,\alpha_2}(s) =}{}\times
 \sum\limits_{l=0}^{n_2}\sum\limits_{j=l}^{l+n_1}\frac{[s]_q^{(j)}[n_1]_q^{(j-l)}[n_2]_q^{(l)}}{[j-l]_q!
[l]_q!}
\left(\frac{-1}{q^{2n_1}\alpha_1}\right)^{j-l}\left(\frac{-1}{q^{2n_2}\alpha_2}\right)^{l}q^{\binom{j-l+1}{2}+\binom{l+1}{2}}
\nonumber
\\
\hphantom{C_{q,n_1,n_2}^{\alpha_1,\alpha_2}(s)}{}
=(-\alpha_1)^{n_1}(-\alpha_2)^{n_2}q^{n_1^2+n_2^2+n_1n_2+\frac{3}{2}\left(\binom{n_1}{2}+\binom{n_2}{2}\right)}
\nonumber
\\
\hphantom{C_{q,n_1,n_2}^{\alpha_1,\alpha_2}(s) =}{}\times
 \sum\limits_{j=0}^{n_1+n_2}\sum\limits_{l=\max (0,j-n_1)}^{\min (j,n_2)}\!\!\frac{[s]_q^{(j)}
[n_1]_q^{(j-l)}[n_2]_q^{(l)}}{[j-l]_q! [l]_q!}\!
\left(\frac{-1}{q^{2n_1}\alpha_1}\right)^{j-l}\!\left(\frac{-1}{q^{2n_2}\alpha_2}\right)^{l}q^{\binom{j-l+1}{2}+\binom{l+1}{2}}
\nonumber
\\
\hphantom{C_{q,n_1,n_2}^{\alpha_1,\alpha_2}(s)  }{}
=
 (-\alpha_1)^{n_1}(-\alpha_2)^{n_2}q^{n_1^2+n_2^2+n_1n_2+\frac{3}{2}\left(\binom{n_1}{2}+\binom{n_2}{2}\right)}
\label{eqq-iden}
\\
\hphantom{C_{q,n_1,n_2}^{\alpha_1,\alpha_2}(s) =}{}\times
\sum\limits_{j=0}^{n_1+n_2}\sum\limits_{l=0}^{j}\frac{[s]_q^{(j)}[n_1]_q^{(j-l)}
[n_2]_q^{(l)}}{[j-l]_q! [l]_q!}
\left(\frac{-1}{q^{2n_1}\alpha_1}\right)^{j-l}\left(\frac{-1}{q^{-2n_2}\alpha_2}\right)^{l}q^{\binom{j-l+1}{2}+\binom{l+1}{2}}.
\nonumber
\end{gather}
This expression is useful for computing some recurrence coef\/f\/icients (see Section~\ref{rr}).

\section[High-order $q$-dif\/ference equation]{High-order $\boldsymbol{q}$-dif\/ference equation}\label{df-eq}

In this section we will f\/ind a~lowering operator for the $q$-Charlier multiple orthogonal polyno\-mials.
Then we will combine it with the raising operators~\eqref{ROpqMCharlier} to get an $(r+1)$-order $q$-dif\/ference equation
(in the same fashion that~\cite{arvesu-esposito} and~\cite{lee}).
\begin{Lemma}
Let $\mathbb{V}$ be the linear subspace of polynomials $Q(s)$ on the lattice $x(s)$ of degree at most $\vert
\vec{n}\vert-1$ defined by the following conditions
\begin{gather*}
\sum\limits_{s=0}^{\infty}Q(s) [s]_{q}^{(k)}\upsilon_{q}^{q\alpha_{j}}(s)\bigtriangledown x_{1}(s)=0,
\qquad
0\leq k\leq n_{j}-2
\qquad
\text{and}
\qquad
j=1,\ldots,r.
\end{gather*}
Then, the system $\{C_{q,\vec{n}-\vec{e}_{i}}^{\vec{\alpha}_{i,q}}(s) \}_{i=1}^{r}$, where
$\vec{\alpha}_{i,q}=(\alpha_{1},\ldots,q\alpha_{i},\ldots,\alpha_{r})$, is a~basis for $\mathbb{V}$.
\label{LI}
\end{Lemma}

\begin{proof}

From orthogonality relations
\begin{gather*}
\sum\limits_{s=0}^{\infty}C_{q,\vec{n}-\vec{e}_{i}}^{\vec{\alpha}_{i,q}}(s)
[s]_{q}^{(k)}\upsilon_{q}^{q\alpha_{j}}(s) \bigtriangledown x_{1}(s) =0,
\qquad
0\leq k\leq n_{j}-2,
\qquad
j=1,\ldots,r,
\end{gather*}
we have that polynomials $C_{q,\vec{n}-\vec{e}_{i}}^{\vec{\alpha}_{i,q}}(s)$, $i=1,\ldots,r$, belong to $\mathbb{V}$.

Now, aimed to get a~contradiction, let us assume that there exists constants $\lambda_{i}$, $i=1,\ldots,r$, such that
\begin{gather*}
\sum\limits_{i=1}^{r}\lambda_{i}C_{q,\vec{n}-\vec{e}_{i}}^{\vec{\alpha}_{i,q}}(s) =0,
\qquad
\text{where}
\qquad
\sum\limits_{i=1}^{r}\abs{\lambda_{i}}>0.
\end{gather*}
Then, multiplying the previous equation by $[s]_{q}^{(n_{k}-1)}\upsilon_{q}^{\alpha_{k}}(s)\bigtriangledown
x_{1}(s)$ and then taking summation on~$s$ from $0$ to~$\infty$, one gets
\begin{gather*}
\sum\limits_{i=1}^{r}\lambda_{i}\sum\limits_{s=0}^{\infty}C_{q,\vec{n}-\vec{e}_{i}}^{\vec{\alpha}_{i,q}}(s)
[s]_{q}^{(n_{k}-1)}\upsilon_{q}^{\alpha_{k}}(s)\bigtriangledown x_{1}(s) =0.
\end{gather*}
Thus, taking into account relations
\begin{gather}
\sum\limits_{s=0}^{\infty}C_{q,\vec{n}-\vec{e}_{i}}^{\vec{\alpha}_{i,q}}(s)
[s]_{q}^{(n_{k}-1)}\upsilon_{q}^{\alpha_{k}}(s) \bigtriangledown x_{1}(s) =c\delta_{i,k},
\qquad
c\in\mathbb{R}\setminus\{0\},
\label{Orthog_Cond}
\end{gather}
we deduce that $\lambda_{k}=0$ for $k=1,\ldots,r$.
Here $\delta_{i,k}$ represents the Kronecker delta symbol.
Therefore, $\{C_{q,\vec{n}-\vec{e}_{i}}^{\vec{\alpha}_{i,q}}(s) \}_{i=1}^{r}$ is linearly independent in
$\mathbb{V}$.
Furthermore, we know that any polynomial of $\mathbb{V}$ can be determined with $\vert \vec{n}\vert$
coef\/f\/icients while $(\vert \vec{n}\vert -r)$ linear conditions are imposed on $\mathbb{V}$,
consequently the dimension of $\mathbb{V}$ is at most~$r$.
Hence, the system $\{C_{q,\vec{n}-\vec{e}_{i}}^{\vec{\alpha}_{i,q}}(s)\}_{i=1}^{r}$ spans $\mathbb{V}$, which
completes the proof.
\end{proof}

Now we will prove that operator~\eqref{lower} is indeed a~lowering operator for the sequence of $q$-Charlier multiple
orthogonal polynomials $C_{q,\vec{n}}^{\vec{\alpha}}(s)$.

\begin{Lemma}
There holds the following relation
\begin{gather}
\Delta C_{q,\vec{n}}^{\vec{\alpha}}(s) =\sum\limits_{i=1}^{r}q^{|\vec{n}|-n_{i}+1/2}
[n_{i}]_{q}^{(1)}C_{q,\vec{n}-\vec{e}_{i}}^{\vec{\alpha}_{i,q}}(s).
\label{Rela_qChar}
\end{gather}
\end{Lemma}

\begin{proof}
Using summation by parts we have
\begin{gather}
\sum\limits_{s=0}^{\infty}\Delta C_{q,\vec{n}}^{\vec{\alpha}}(s)[s]_{q}^{(k)}\upsilon_{q}^{q\alpha_{j}}(s) \bigtriangledown x_{1}(s)
=-\sum\limits_{s=0}^{\infty}C_{q,\vec{n}}^{\vec{\alpha}}(s)\nabla\big[[s]_{q}^{(k)}\upsilon_{q}^{q\alpha_{j}}(s)\big] \bigtriangledown x_{1}(s)
\nonumber
\\
\hphantom{\sum\limits_{s=0}^{\infty}\Delta C_{q,\vec{n}}^{\vec{\alpha}}(s)[s]_{q}^{(k)}\upsilon_{q}^{q\alpha_{j}}(s) \bigtriangledown x_{1}(s)}{}
=-\sum\limits_{s=0}^{\infty}C_{q,\vec{n}}^{\vec{\alpha}}(s)\varphi_{j,k}(s)\upsilon_{q}^{\alpha_{j}}(s) \bigtriangledown x_{1}(s),
\label{inte-1}
\end{gather}
where
\begin{gather*}
\varphi_{j,k}(s) =q^{1/2}[s]_{q}^{(k)}-q^{-1/2}\frac{x(s)}{\alpha_{j}}[s-1]_{q}^{(k)},
\end{gather*}
is a~polynomial of degree $\leq k+1$ in the variable $x(s)$.
Consequently, from the orthogonality conditions~\eqref{NOrthC} we get
\begin{gather*}
\sum\limits_{s=0}^{\infty}\Delta C_{q,\vec{n}}^{\vec{\alpha}}(s) [s]_{q}^{(k)}\upsilon_{q}^{q\alpha_{j}}(s)
\bigtriangledown x_{1}(s)=0,
\qquad
0\leq k\leq n_{j}-2,
\qquad
j=1,\ldots,r.
\end{gather*}
Hence, from Lemma~\ref{LI}, $\Delta C_{q,\vec{n}}^{\vec{\alpha}}(s) \in \mathbb{V}$.
Moreover, $\Delta C_{q,\vec{n}}^{\vec{\alpha}}(s)$ can univocally be expressed as a~linear combination of polynomials
$\{C_{q,\vec{n}-\vec{e}_{i}}^{\vec{\alpha}_{i,q}}(s) \}_{i=1}^{r}$, i.e.~
\begin{gather}
\Delta C_{q,\vec{n}}^{\vec{\alpha}}(s) =\sum\limits_{i=1}^{r}\beta_{i}C_{q,\vec{n}-\vec{e}_{i}}^{\vec{\alpha}_{i,q}}(s),
\qquad
\sum\limits_{i=1}^{r}\abs{\beta_{i}}>0.
\label{eq-del}
\end{gather}
Multiplying both sides of the equation~\eqref{eq-del} by $[s]_{q}^{(n_{k}-1)}\upsilon_{q}^{q\alpha_{k}}(s)
\bigtriangledown x_{1}(s)$ and using relations~\eqref{Orthog_Cond} one has
\begin{gather}
\sum\limits_{s=0}^{\infty}\Delta C_{q,\vec{n}}^{\vec{\alpha}}(s)[s]_{q}^{(n_{k}-1)}\upsilon_{q}^{q\alpha_{k}}(s)\bigtriangledown x_{1}(s)
=\sum\limits_{i=1}^{r}\beta_{i}\sum\limits_{s=0}^{\infty}C_{q,\vec{n}-\vec{e}_{i}}^{\vec{\alpha}_{i,q}}(s)
[s]_{q}^{(n_{k}-1)}\upsilon_{q}^{q\alpha_{k}}(s) \bigtriangledown x_{1}(s)
\nonumber
\\
\hphantom{\sum\limits_{s=0}^{\infty}\Delta C_{q,\vec{n}}^{\vec{\alpha}}(s)[s]_{q}^{(n_{k}-1)}\upsilon_{q}^{q\alpha_{k}}(s)\bigtriangledown x_{1}(s)}{}
=\beta_{k}\sum\limits_{s=0}^{\infty}C_{q,\vec{n}-\vec{e}_{k}}^{\vec{\alpha}_{i,q}}(s)[s]_{q}^{(n_{k}-1)}\upsilon_{q}^{q\alpha_{k}}(s)
\bigtriangledown x_{1}(s).
\label{Ident_I}
\end{gather}
If we replace $[s]_{q}^{(k)}$ by $[s]_{q}^{(n_k-1)}$ in the left-hand side of equation~\eqref{inte-1}, then left-hand
side of equation~\eqref{Ident_I} transforms into relation
\begin{gather}
\sum\limits_{s=0}^{\infty}\Delta C_{q,\vec{n}}^{\vec{\alpha}}(s)[s]_{q}^{(n_{k}-1)}\upsilon_{q}^{q\alpha_{k}}(s) \bigtriangledown x_{1}(s)
=-\sum\limits_{s=0}^{\infty}C_{q,\vec{n}}^{\vec{\alpha}}(s)\varphi_{k,n_{k}-1}(s)\upsilon_{q}^{\alpha_{k}}(s) \bigtriangledown x_{1}(s)
\nonumber
\\
\hphantom{\sum\limits_{s=0}^{\infty}\Delta C_{q,\vec{n}}^{\vec{\alpha}}(s)[s]_{q}^{(n_{k}-1)}\upsilon_{q}^{q\alpha_{k}}(s) \bigtriangledown x_{1}(s)}{}
=\frac{q^{-1/2}}{\alpha_{k}}\sum\limits_{s=0}^{\infty}C_{q,\vec{n}}^{\vec{\alpha}}(s)[s]_{q}^{(n_{k})}\upsilon_{q}^{\alpha_{k}}(s)
\bigtriangledown x_{1}(s).
\label{eqcha-ra}
\end{gather}
Here we have used that $x(s)[s-1]_{q}^{(n_{k}-1)}=[s]_{q}^{(n_{k})}$ to get
$\varphi_{k,n_{k}-1}(s) =- (q^{-1/2}/\alpha_{k})[s]_{q}^{(n_{k})}+ \text{lower terms}$.

On the other hand, from~\eqref{ROpqMCharlier} one has that
\begin{gather}
\frac{1}{\alpha_{k}}\upsilon_{q}^{\alpha_{k}}(s) C_{q,\vec{n}}^{\vec{\alpha}}(s) =-q^{|\vec{n}|-1/2}\nabla
\big[\upsilon_{q}^{q\alpha_{k}}(s) C_{q,\vec{n}-\vec{e}_{k}}^{\vec{\alpha}_{i,q}}(s)\big].
\label{eqcha-ra1}
\end{gather}
Then, by conveniently substituting~\eqref{eqcha-ra1} in the right-hand side of equation~\eqref{eqcha-ra} and using once
more summation by parts, we get
\begin{gather*}
\begin{split}
& \sum\limits_{s=0}^{\infty}\Delta C_{q,\vec{n}}^{\vec{\alpha}}(s)
[s]_{q}^{(n_{k}-1)}\upsilon_{q}^{q\alpha_{k}}(s) \bigtriangledown x_{1}(s)
 =-q^{|\vec{n}|-1}\sum\limits_{s=0}^{\infty}[s]_{q}^{(n_{k})}\nabla \big[\upsilon_{q}^{q\alpha_{k}}(s)
C_{q,\vec{n}-\vec{e}_{k}}^{\vec{\alpha}_{i,q}}(s) \big] \bigtriangledown x_{1}(s)
\\
& \hphantom{\sum\limits_{s=0}^{\infty}\Delta C_{q,\vec{n}}^{\vec{\alpha}}(s)
[s]_{q}^{(n_{k}-1)}\upsilon_{q}^{q\alpha_{k}}(s) \bigtriangledown x_{1}(s)}{}
 =q^{|\vec{n}|-1}\sum\limits_{s=0}^{\infty}C_{q,\vec{n}-\vec{e}_{k}}^{\vec{\alpha}_{i,q}}(s) \Delta
\big[[s]_{q}^{(n_{k})}\big] \upsilon_{q}^{q\alpha_{k}}(s) \bigtriangledown x_{1}(s).
\end{split}
\end{gather*}
Since $\Delta [s]_{q}^{(n_{k})}=q^{3/2-n_{k}}[n_{k}]_{q}^{(1)}[s]_{q}^{(n_{k}-1)}$ we
f\/inally have
\begin{gather*}
\sum\limits_{s=0}^{\infty}\Delta C_{q,\vec{n}}^{\vec{\alpha}}(s)
[s]_{q}^{(n_{k}-1)}\upsilon_{q}^{q\alpha_{k}}(s) \bigtriangledown x_{1}(s)
\\
\qquad{}
=q^{|\vec{n}|-n_{k}+1/2}[n_{k}]_{q}^{(1)}\sum\limits_{s=0}^{\infty}C_{q,\vec{n}-\vec{e}_{k}}^{\vec{\alpha}_{i,q}}(s)
[s]_{q}^{(n_{k}-1)}\upsilon_{q}^{q\alpha_{k}}(s) \bigtriangledown x_{1}(s).
\end{gather*}
Therefore, comparing this equation with~\eqref{Ident_I} we obtain the coef\/f\/icients in the expansion~\eqref{eq-del},
i.e.~
\begin{gather*}
\beta_{k}=q^{|\vec{n}|-n_{k}+1/2}[n_{k}]_{q}^{(1)},
\end{gather*}
which proves relation~\eqref{Rela_qChar}.
\end{proof}

\begin{Theorem}
The $q$-Charlier multiple orthogonal polynomial $C_{q,\vec{n}}^{\vec{\alpha}}(s)$ satisfies the following
$(r+1)$-order $q$-difference equation
\begin{gather}
\prod\limits_{i=1}^{r}\mathcal{D}_{q}^{\alpha_{i}}\Delta C_{q,\vec{n}}^{\vec{\alpha}}(s)
=-\sum\limits_{i=1}^{r}q^{\vert \vec{n}\vert -n_{i}+1}[n_{i}]_{q}^{(1)}
\prod\limits_{\substack{j=1
\\
j\neq i}}^{r}\mathcal{D}_{q}^{\alpha_{j}}C_{q,\vec{n}}^{\vec{\alpha}}(s).
\label{q-DEquation}
\end{gather}
\end{Theorem}

\begin{proof}
Since operators~\eqref{ROpqMCharlier} are commuting, we write
\begin{gather}
\prod\limits_{i=1}^{r}\mathcal{D}_{q}^{\alpha_{i}}
=\left(\prod\limits_{\substack{j=1\\j\neq i}}^{r} \mathcal{D}_{q}^{\alpha_{j}}\right) \mathcal{D}_{q}^{\alpha_{i}},
\label{interm}
\end{gather}
and then using~\eqref{ROpqMCharlier}, by acting on equation~\eqref{Rela_qChar} with the product of
operators~\eqref{interm}, we obtain the following relation
\begin{gather*}
\prod\limits_{i=1}^{r}\mathcal{D}_{q}^{\alpha_{i}}\Delta C_{q,\vec{n}}^{\vec{\alpha}}(s)
 = \sum\limits_{i=1}^{r}q^{|\vec{n}|-n_{i}+1/2}[n_{i}]_{q}^{(1)}\left(\prod\limits_{\substack{j=1\\j\neq i}}^{r}\mathcal{D}_{q}^{\alpha_{j}}\right)
\mathcal{D}_{q}^{\alpha_{i}}C_{q,\vec{n}-\vec{e}_{i}}^{\vec{\alpha}_{i,q}}(s)
\\
\hphantom{\prod\limits_{i=1}^{r}\mathcal{D}_{q}^{\alpha_{i}}\Delta C_{q,\vec{n}}^{\vec{\alpha}}(s)}{}
 = -\sum\limits_{i=1}^{r}q^{\vert \vec{n}\vert -n_{i}+1}[n_{i}
]_{q}^{(1)}\prod\limits_{\substack{j=1
\\
j\neq i}}^{r} \mathcal{D}_{q}^{\alpha_{j}}C_{q,\vec{n}}^{\vec{\alpha}}(s),
\end{gather*}
which proves~\eqref{q-DEquation}.
\end{proof}

\section{Recurrence relation}\label{rr}

In this section we will study two types of recurrence relations, namely the nearest neighbor recurrence relation for
any multi-index $\vec{n}$ and a~step-line recurrence relation for $\vec{n}=(n_1,n_2)$.

For the main result of this section we will make use of the following lemma.

\begin{Lemma}\label{Coro-q-Ch}
Let $n_{i}$ be a~positive integer and let $f(s)$ be a~function defined on the discrete variable~$s$.
The following relation is valid
\begin{gather}
\mathcal{C}_{q,n_{i}}^{\alpha_{i}}(x(s) f(s)) =q^{-n_{i}+1/2}x(n_{i})
(\alpha_{i})^{-s}\nabla^{n_{i}-1} (\alpha_{i}q^{n_{i}})^{s}f(s)\! +q^{-n_{i}}[x(s)\!
-x(n_{i}) ] \mathcal{C}_{q,n_{i}}^{\alpha_{i}}f(s),\!\!\!
\label{eq-lem-p}
\end{gather}
where difference operator $\mathcal{C}_{q,n_{i}}^{\alpha_{i}}$ is given in~\eqref{Op-q-Ch}.
\end{Lemma}

\begin{proof}
Let us act $n_{i}$-times with backward dif\/ference operators~\eqref{nbackward} on the product of functions $x(s)f(s)$.
Assume momentarily that $n_i\geq N>1$,
\begin{gather}
\nabla^{n_{i}}x(s) f(s)=\nabla^{n_{i}-1} [\nabla x(s) f(s) ] =\nabla^{n_{i}-1}\big[q^{-1/2}f(s)+x(s-1) \nabla f(s) \big]
\nonumber\\
\phantom{\nabla^{n_{i}}x(s) f(s)}{}
=q^{-1/2}\nabla^{n_{i}-1}f(s) +\nabla^{n_{i}-1}[x(s-1) \nabla f(s) ]
\nonumber\\
\phantom{\nabla^{n_{i}}x(s) f(s)}
=q^{-1/2}\nabla^{n_{i}-1}f(s) +\nabla^{n_{i}-2}[\nabla x(s-1) \nabla f(s) ].
\label{rety}
\end{gather}
Repeating this process~-- but on the second term of the right-hand side of equation~\eqref{rety}~-- yields
\begin{gather*}
\nabla^{n_{i}}x(s) f(s)= \big[q^{1/2-n_{i}}+\cdots +q^{-5/2}+q^{-3/2}+q^{-1/2}\big] \nabla^{n_{i}-1}f(s)
+x(s-n_{i}) \nabla^{n_{i}}f(s)
\\
\phantom{\nabla^{n_{i}}x(s) f(s)}{}
 = q^{1/2-n_{i}}x(n_i) \nabla^{n_{i}-1}f(s) +x(s-n_{i}) \nabla^{n_{i}}f(s),
\end{gather*}
or equivalently,
\begin{gather}
\nabla^{n_{i}}x(s) f(s) =q^{-n_{i}+1/2}x(n_{i}) \nabla^{n_{i}-1}f(s) +q^{-n_{i}}[x(s)
-x(n_{i}) ] \nabla^{n_{i}}f(s),
\qquad
n_{i}\geq 1.
\label{eq-semi}
\end{gather}
Now, aimed to involve dif\/ference operator $\mathcal{C}_{q,n_{i}}^{\alpha_{i}}$ in the above equation, we multiply
equation~\eqref{eq-semi} from the left side by $(\alpha_{i})^{-s}$ and replace $f(s)$ by
$(\alpha_{i}q^{n_{i}})^{s}f(s)$.
Thus, equation~\eqref{eq-semi} transforms into~\eqref{eq-lem-p}, which completes the proof.
\end{proof}

\begin{Theorem}
The $q$-Charlier multiple orthogonal polynomials satisfy the following $(r+2)$-term recurrence relation
\begin{gather}
x(s) C_{q,\vec{n}}^{\vec{\alpha}}(s) =C_{q,\vec{n}+\vec{e}_{k}}^{\vec{\alpha}}(s)
+\!\left(\sum\limits_{i=1}^{r}q^{\vert
\vec{n}\vert_{i}}x(n_{i})\bigg[(q-1)\bigg(\alpha_{i}q^{\sum\limits_{j=i}^{r}n_{j}}\bigg)\!+1\bigg]\!
+\alpha_{k}q^{|\vec{n}|+n_k+1}\right)\! C_{q,\vec{n}}^{\vec{\alpha}}(s)
\nonumber
\\
\hphantom{x(s) C_{q,\vec{n}}^{\vec{\alpha}}(s) =}{}
+\sum\limits_{i=1}^{r}q^{\vert
\vec{n}\vert_{i}}x(n_{i})\bigg[(q-1)\bigg(\alpha_{i}q^{\sum\limits_{j=i}^{r}n_{j}}\bigg)+1\bigg]\alpha_iq^{|\vec{n}|+n_i-1}
C_{q,\vec{n}-\vec{e}_{i}}^{\vec{\alpha}}(s),
\label{q-Rrelation}
\end{gather}
where $\vert \vec{n}\vert_{i}=n_{1}+\cdots +n_{i-1}$, $\vert \vec{n}\vert_{1}=0$.
\end{Theorem}

\begin{proof}
We will use Lemma~\ref{Coro-q-Ch}.
First, let us consider equation
\begin{gather*}
(\alpha_{k})^{-s}\nabla^{n_{k}+1}(\alpha_{k}q^{n_k+1})^{s} \frac{1}{\Gamma_{q}(s+1)}
=(\alpha_{k})^{-s}\nabla^{n_{k}} \left[q^{-s+1/2}\bigtriangledown
\frac{(\alpha_{k}q^{n_k+1})^{s}}{\Gamma_{q}(s+1)}\right]
\\
\hphantom{(\alpha_{k})^{-s}\nabla^{n_{k}+1}(\alpha_{k}q^{n_k+1})^{s} \frac{1}{\Gamma_{q}(s+1)}}{}
=q^{1/2}(\alpha_{k})^{-s}\nabla^{n_{k}}(\alpha_{k}q^{n_k})^{s} \frac{1}{\Gamma_{q}(s+1)}
\\
\hphantom{(\alpha_{k})^{-s}\nabla^{n_{k}+1}(\alpha_{k}q^{n_k+1})^{s} \frac{1}{\Gamma_{q}(s+1)}=}{}
-(\alpha_{k}q^{n_k+1})^{-1}q^{1/2}(\alpha_{k})^{-s}
\nabla^{n_{k}}(\alpha_{k}q^{n_k})^{s} x(s)\frac{1}{\Gamma_{q}(s+1)},
\end{gather*}
which can be rewritten in terms of dif\/ference operators~\eqref{Op-q-Ch} as follows
\begin{gather}
\mathcal{C}_{q,n_{k}+1}^{\alpha_{k}}\frac{1}{\Gamma_{q}(s+1)}=q^{1/2}
\mathcal{C}_{q,n_{k}}^{\alpha_{k}}\frac{1}{\Gamma_{q}(s+1)}-(\alpha_{k}q^{n_k+1})^{-1}q^{1/2}
\mathcal{C}_{q,n_{k}}^{\alpha_{k}}x(s) \frac{1}{\Gamma_{q}(s+1)}.
\label{eqsvec}
\end{gather}
Since operators~\eqref{Op-q-Ch} are commuting the multiplication of equation~\eqref{eqsvec} from the left-hand side~by
the product $\left(\prod\limits_{\substack{j=1\\j\neq k}}^{r}\mathcal{C}_{q,n_{j}}^{\alpha_{j}}\right)$ yields
\begin{gather}
\mathcal{C}_{q,\vec{n}}^{\vec{\alpha}}x(s)
\frac{1}{\Gamma_{q}(s+1)}=\alpha_{k}q^{n_k+1}\big(\mathcal{C}_{q,\vec{n}}^{\vec{\alpha}}
-q^{-1/2}\mathcal{C}_{q,\vec{n} +\vec{e}_{k}}^{\vec{\alpha}}\big)\frac{1}{\Gamma_{q}(s+1)}.
\label{eq-rr-interm}
\end{gather}

Second, let us recursively use expression~\eqref{eq-lem-p} involving the product of~$r$ dif\/ference operators
$\mathcal{C}_{q,n_{1}}^{\alpha_{1}},\ldots,\mathcal{C}_{q,n_{r}}^{\alpha_{r}}$ acting on the function $f(s)
=1/\Gamma_{q}(s+1)$.
Thus, we have
\begin{gather}
q^{\vert \vec{n}\vert}\mathcal{C}_{q,\vec{n}}^{\vec{\alpha}} x(s)\frac{1}{\Gamma_{q}(s+1)}
\nonumber\\
\qquad
=q^{1/2}\sum\limits_{i=1}^{r}q^{\vert \vec{n}\vert_{i}}x(n_{i})
\bigg[(q-1)\bigg(\alpha_{i}q^{\sum\limits_{j=i}^{r}n_{j}}\bigg)+1\bigg]
\prod\limits_{j=1}^{r}\mathcal{C}_{q,n_{j}-\delta_{j,i}}^{\alpha_{j}}\frac{1}{\Gamma_{q}(s+1)}
\nonumber\\
\qquad\hphantom{=}{}
+\left(x(s) -\sum\limits_{i=1}^{r}q^{\vert \vec{n}
\vert_{i}}x(n_{i})\bigg[(q-1)\bigg(\alpha_{i}q^{\sum\limits_{j=i}^{r}n_{j}}\bigg)+1\bigg]
\right) \mathcal{C}_{q,\vec{n}}^{\vec{\alpha}}\frac{1}{\Gamma_{q}(s+1)}.
\label{eq-rr-interm1}
\end{gather}
Hence, by using expressions~\eqref{eq-rr-interm}, \eqref{eq-rr-interm1} one gets
\begin{gather*}
x(s)\mathcal{C}_{q,\vec{n}}^{\vec{\alpha}}\frac{1}{\Gamma_{q}(s+1)}
=-\alpha_{k}q^{|\vec{n}|+n_k+1}q^{-1/2}\mathcal{C}_{q,\vec{n}+\vec{e}_{k}}^{\vec{\alpha}}
\frac{1}{\Gamma_{q}(s+1)}
\\
\qquad{}
+\left(\sum\limits_{i=1}^{r}q^{\vert\vec{n}\vert_{i}}x(n_{i})\bigg[(q-1)\bigg(\alpha_{i}q^{\sum\limits_{j=i}^{r}n_{j}}\bigg)+1\bigg]
+\alpha_{k}q^{|\vec{n}|+n_k+1}\right) \mathcal{C}_{q,\vec{n}}^{\vec{\alpha}}\frac{1}{\Gamma_{q}(s+1)}
\\
\qquad{}
-q^{1/2}\sum\limits_{i=1}^{r}q^{\vert\vec{n}\vert_{i}}x(n_{i})
\bigg[(q-1)\bigg(\alpha_{i}q^{\sum\limits_{j=i}^{r}n_{j}}\bigg)+1\bigg]\prod\limits_{j=1}^{r}\mathcal{C}_{q,n_{j}
-\delta_{j,i}}^{\alpha_{j}}\frac{1}{\Gamma_{q}(s+1)}.
\end{gather*}

Finally, multiplying from the left both sides of the previous expression~by
$\mathcal{K}_{q}^{\vec{n},\vec{\alpha}}\Gamma_{q}(s+1)$ and using Rodrigues-type formula~\eqref{RFormula} we
obtain~\eqref{q-Rrelation}, which completes the proof.
\end{proof}

Observe that other recurrence relations dif\/ferent from the above nearest neighbor recurrence
relation~\eqref{q-Rrelation} can be obtained.
Indeed, from~\eqref{eqq-iden} a~$4$-term recurrence relation for $\vec{n}=(n_1, n_2)$ can be obtained.
In~\cite{arvesu_vanAssche} it was given an approach for the recurrence relations of some discrete multiple orthogonal
polynomials.
This approach make use of the Rodrigues-type formulas along with some series representations (Kamp\'e de F\'eriet
series) for multiple orthogonal polynomials.
Here we proceed in the same fashion.

Considering the expansion
\begin{gather*}
C_{q,n_1,n_2}^{\alpha_1,\alpha_2}(s)=\sum\limits_{j=0}^{n_1+n_2}c_{q,n_1,n_2}^{(j)}[s]_q^{(j)},
\end{gather*}
the coef\/f\/icients $c_{q,n_1,n_2}^{(j)}$ can be used to compute the recurrence coef\/f\/icients in
\begin{gather*}
\begin{split}
& x(s)P_{q,n_1,n_2}(s)=q^{n_1+n_2}P_{q,n_1,n_2+1}(s)+b_{q,n_1,n_2}P_{q,n_1,n_2}(s)+c_{q,n_1,n_2}P_{q,n_1,n_2-1}(s)
\\
& \phantom{x(s)P_{q,n_1,n_2}(s)=}{}
+d_{q,n_1,n_2}P_{q,n_1-1,n_2-1}(s),
\end{split}
\end{gather*}
where $P_{q,n_1,n_2}(s)=C_{q,n_1,n_2}^{\alpha_1,\alpha_2}(s)$.
Indeed,
\begin{gather*}
b_{q,n_1,n_2} = q^{n_1+n_2}\big(q^{-1}c_{q,n_1,n_2}^{(n_1+n_2-1)}-c_{q,n_1,n_2+1}^{(n_1+n_2)}\big)+x(n_1+n_2),
\\
c_{q,n_1,n_2}=q^{n_1+n_2}\big(q^{-2}c_{q,n_1,n_2}^{(n_1+n_2-2)}+q^{-(n_1+n_2)}c_{q,n_1,n_2}^{(n_1+n_2-1)}x
(n_1+n_2-1)  \\
\hphantom{c_{q,n_1,n_2}=}{}
- q^{-(n_1+n_2)}b_{q,n_1,n_2}c_{q,n_1,n_2}^{(n_1+n_2-1)}-c_{q,n_1,n_2+1}^{(n_1+n_2-1)}\big),
\notag\\
d_{q,n_1,n_2} = q^{n_1+n_2}\big(q^{-3}c_{q,n_1,n_2}^{(n_1+n_2-3)}+q^{-(n_1+n_2)}c_{q,n_1,n_2}^{(n_1+n_2-2)}x(n_1+n_2-2) \\
\hphantom{d_{q,n_1,n_2} =}{}
- q^{-(n_1+n_2)}c_{q,n_1,n_2}c_{q,n_1,n_2-1}^{(n_1+n_2-2)}-q^{-(n_1+n_2)}b_{q,n_1,n_2}c_{q,n_1,n_2}^{(n_1+n_2-2)}
-c_{q,n_1,n_2+1}^{(n_1+n_2-2)}\big).
\end{gather*}

From the explicit expression~\eqref{eqq-iden} we then get, after some calculations, that when $q\rightarrow 1$ we
recover the recurrence coef\/f\/icients given in~\cite{arvesu_vanAssche}.

\section{Conclusions}\label{conclu}

In summary, let us recall some of our results.
We have introduced a~new family of special functions, i.e.\ $q$-Charlier multiple orthogonal polynomials.
For these polynomials we have obtained a~Rodrigues-type formula~\eqref{RFormula} as well as their explicit
representation in terms of a~$q$-analogue of the second of Appell's hypergeometric functions of two
variables~\eqref{phi2_charli}.
Moreover, these polynomials are common eigenfunctions of two dif\/ferent $(r+1)$-order dif\/ference operators,
namely~\eqref{q-DEquation} and~\eqref{q-Rrelation}.
In the limiting situation $q\rightarrow 1$ we recover the corresponding structural relations for multiple Charlier
polynomials~\cite{arvesu_vanAssche}.
Indeed, our relations~\eqref{RFormula},~\eqref{q-DEquation}, and~\eqref{q-Rrelation} transform
into~\eqref{Rodrigues-multi},~\eqref{opdi-1}, and~\eqref{AR}, respectively.

Our algebraic approach for the nearest neighbor recurrence relation~\eqref{q-Rrelation} does not require to introduce
type I multiple orthogonality~\cite{Assche_neighbor}.
Indeed, the $q$-dif\/ference operators involved in the Rodrigues-type formula constitute the key-ingredient in our
approach.

A description of the main term of the logarithm asymptotics of $q$-Charlier multiple orthogonal polynomials deserves our
future attention.
We expect to give it in terms of an algebraic function (see~\cite{apt-arv}).
These results yield the Cauchy transform of the weak-star limit for scaling zero counting measure of the polynomials.
In addition, the zero distribution of these type II multiple orthogonal polynomials will be studied.

\subsection*{Acknowledgements}
The research of J.~Arves\'u was partially supported by the research grant MTM2012-36732-C03-01 (Ministerio de
Econom\'{\i}a y Competitividad) of Spain.

\pdfbookmark[1]{References}{ref}
\LastPageEnding

\end{document}